\documentclass[12pt,a4paper, english]{paper}  
\usepackage{babel}
\usepackage[algoruled,vlined,linesnumbered]{algorithm2e}
\usepackage{graphicx}%
\usepackage{multirow}%
\usepackage{amsmath,amssymb,amsfonts}%
\usepackage{amsthm}%

\newtheorem{theorem}{Theorem}[section]

\newtheorem{definition}[theorem]{Definition}

\newtheorem{corollary}[theorem]{Corollary}
\newtheorem{lemma}[theorem]{Lemma}

\usepackage[normalem]{ulem}

\usepackage[arc,all,color]{xy}
\newenvironment{grafo}[1]{\begin{xy}0;<#1,0cm>: }{\end{xy}}
\newcommand{\NodoO}[5][5mm]{\POS *=<#1>[o]{#3}*\frm{o}="#2"+#4(1.5)*{#5}}
\newcommand{\ArcoD}[3][]{\ar@{->}"#2";"#3"#1}
\newcommand{\ArcoO}[3][]{\ar@/^/@{->}"#2";"#3"#1}
\newcommand{\ArcoA}[3][]{\ar@/_/@{->}"#2";"#3"#1}
\def\uno{\mathbf 1}
\newcommand{\G}{\mathcal{G}}
\newcommand\R{\mathbb{R}}

\renewcommand{\vec}{\mathrm{vec}}

\newcommand\wh{\widehat}

\begin{document}

\title{Role extraction by matrix equations
and generalized random walks}
%%% \footnote{Last update: \today}}

\author{Dario Fasino}
%ORCID: 0000-0001-7682-0660

\maketitle

\abstract{
The nodes in a network can be grouped into 'roles' based on similar connection patterns. This is usually achieved by defining a pairwise node similarity matrix and then clustering rows and columns of this matrix.
%The nodes in a network may be grouped into equivalence classes according to the role they play, depending on their connections with nodes in either the same role or different roles. To this goal, one firstly defines a similarity 
%matrix whose entries quantify some kind of similarity measure between node pairs. Subsequently, the network roles are 
%detected by clustering rows and columns of this matrix.
This paper presents a new
similarity matrix for solving role extraction problems
in directed networks, which is defined as the solution of a 
matrix equation and computes node similarities based on random walks that can proceed both along the link direction 
and in the opposite direction.
The resulting node similarity measure shows remarkable performance in role extraction tasks
on directed networks with heterogeneous node degree distributions.


}

%\keywords{Role extraction, similarity measure, matrix equation, complex networks, block modeling, random walks}
%
%
%\pacs[MSC Classification]{05C50, 05C85, 65F45, 68R10}

% -------------------------------------------------------
\section{Introduction}

Consider the following graph and assume that it represents a trading network: The nodes correspond to trading companies and the links describe the movement of some goods between them.
\begin{equation}   \label{eq:6nodes}
\begin{grafo}{4.5mm}
   (0,-1.2)\NodoO[3mm]{1}{}{C}{_1},
   (2,-.8)\NodoO[3mm]{2}{}{C}{_2},
   (4,-1.2)\NodoO[3mm]{3}{}{C}{_3},
   (0,1.3)\NodoO[3mm]{4}{}{C}{_4},
   (2,.9)\NodoO[3mm]{5}{}{C}{_5},
   (4,1.3)\NodoO[3mm]{6}{}{C}{_6},   
   \ArcoD[]{1}{2}
   \ArcoO[]{1}{5}
   \ArcoD[]{4}{5}
   \ArcoD[]{2}{3}
   \ArcoD[]{2}{6}
   \ArcoD[]{5}{6}
   \ArcoA[]{2}{5}
   \ArcoA[]{5}{2}
   \end{grafo}
\end{equation} 
Looking at the connections, it is reasonable to 
argue that nodes 1 and 4 represent producers, nodes 
2 and 5 are brokers or dealers, while nodes 3 and 6 
correspond to consumers. 
This simple example, borrowed from \cite{ReiWhi}, illustrates 
the task of role extraction: We are given a directed 
graph or network
and we want to partition its nodes into similarity classes,
called roles,
on the basis of their connection patterns, in a way to 
group the nodes that share the same role or 
behavior within the network. 
As shown by the simple example above, determining the similarity between vertices, or identifying other vertices with similar characteristics to particular vertices, is useful for understanding not only the meso-scale structure but also the functioning of complex networks.

\medskip
Role extraction, also known as block modelling \cite{SBMeasy,ReiWhi,WasFau}, 
aims to extract a simplified, coarse-grained representation of a large, complex network. 
As in the community detection problem \cite{Pan+10,Communities}, 
the purpose of role extraction is to partition nodes in a network into coherent groups.
However, 
the problem is approached from different perspectives. Indeed, 
community detection focuses on identifying cohesive clusters,
whereas 
role extraction aims to identify structural relationships between different groups of nodes based on the different `densities' of connections within and between the groups. 
Besides communities, the identification of other mesoscale structures, such as quasi-bipartite \cite{Bipartite,FT_LAA18} and core-periphery \cite{CP,FR_Sy20,TH-cp} structures in directed and undirected networks, can also be regarded as special instances of the role extraction problem.

\medskip
Role extraction also has its roots in the sociological literature, as is often the case with network science concepts.
In his 1957 work \cite{Nadel}, S.F.\ Nadel pioneered the systematic study of individual roles and statuses within social networks.
Nowadays, it is common for textbooks on social network analysis to include a section on role extraction or block modelling,
see e.g., \cite{ASNbook24,HanRid,WasFau}.

\medskip
One of the most common approach to the role extraction problem is to start by defining a similarity criterion between pairs of nodes. This is done by constructing a similarity matrix whose entries 
quantify some kind of similarity measure of all node pairs. 
Node similarity can be defined in a number of ways, depending on the context. For example, web pages can be considered similar if they contain similar keywords or links to other similar pages, while individuals in a social network can be considered similar based on shared interests, preferences or acquaintances. In any case, similarity measures take on quantitative terms that define a similarity matrix.
In a later step, roles are derived by grouping rows and columns of the similarity matrix.
This last passage is usually performed 
by applying 
a clustering algorithm to the rows or columns of the similarity matrix or a low-rank approximation of it, see e.g.,
\cite{Bar+13,BroVDo}, typical choices being the $k$-means algorithm
\cite{kmeans} or a Louvain-type community detection algorithm \cite{Louvain},
due to their simplicity, fast convergence 
and solid mathematical foundation.

\medskip
This paper introduces a variant of the similarity matrix originally proposed by A.\ Browet and P.\ Van Dooren
in \cite{BroVDo}, which has been subsequently investigated 
also by other authors \cite{Barbarino22,Marchand+}. Similar to the Browet-Van Dooren matrix, our approach involves solving a matrix equation for which we provide a globally convergent iterative method. Our results suggest that 
role extraction tasks based on this refined similarity matrix achieve superior performance in networks where node degrees exhibit significant heterogeneity.
The next section contains the main definitions and 
tools needed in the rest of the paper, and 
provides a short overview
on the literature on node similarity measures.
Section \ref{sec:GWP} introduces 
Browet-Van Dooren's similarity matrix 
and lies the basis for the main results
in this work, which are collected in Section \ref{sec:main}.
One special property of our approach is the subject 
of Section \ref{sec:SBM}, while Section \ref{sec:numerical}
illustrates our findings by numerical examples.

% -------------------------------------------------------
\section{Background and definitions}  \label{sec:2}

In this section, we introduce some notation and basic concepts that we will extensively use throughout this paper. 
We first consider some general definitions and properties of graph theory \cite{Biggs}, matrix analysis and linear algebra \cite{HorJoh,MeyerBook}. 
A directed graph or digraph $\G = (V,E)$ is defined by a set $V$ of nodes or vertices
and a relation $E\subseteq V\times V$.
For notational simplicity, we assume that 
$V = \{1,\ldots,n\}$. The elements of $E$ are called arcs, links or directed edges.
For simplicity, we write
$i\to j$ in place of the more rigorous formula $(i,j)\in E$,
and say that $(i,j)\in E$ is a link from $i$ to $j$.
When the relation $E$ is symmetric then $\G$ is called 
(undirected) graph.
Moreover, we will use the term network with some liberality, as a synonym for graph or digraph.
The adjacency matrix $A$ of $\G = (V,E)$ is the $n\times n$ matrix such that $A_{ij} = 1$ iff $(i,j)\in E$, $0$ otherwise. 
A weighted digraph is a graph where each link is 
associated with a positive number called weight.
In this case, if $i\to j$ then $A_{ij}$ is the weight of that link.
Note that $A^T$
is the adjacency matrix of the digraph obtained by reversing the links of $\G$.

\medskip
The in-degree of a node $i\in V$ is the number
$d^{in}_i = \sum_{j=1}^n A_{ji}$ and the out-degree 
is given by $d^{out}_i = \sum_{j=1}^n A_{ij}$. 
If $\G$ is not weighted, they coincide with
the number of links incoming to and outgoing from node $i$,
respectively. We will make use of the diagonal 
matrices $D_{in} = \mathrm{Diag}(d^{in}_1,\ldots,d^{in}_n)$ and
$D_{out} = \mathrm{Diag}(d^{out}_1,\ldots,d^{out}_n)$.
The in-neighborhood of $i\in V$
is the set 
$\Gamma^{in}_i = \{j\in V: j\to i\}$.
Analogously, the out-neighborhood of $i\in V$
is the set 
$\Gamma^{out}_i = \{j\in V: i\to j\}$.
If $j\in \Gamma^{in}_i$ then we say that $j$ is a parent or ascendant of $i$,
while if $j\in \Gamma^{out}_i$ then $j$ is a child
or descendant of $i$.
A sequence of nodes $i_0,i_1,\ldots,i_\ell$ 
such that $(i_{k-1},i_k)\in E$
for $k = 1,\ldots, \ell$ is a walk 
of length $\ell$.
We say that $i_0$ is the source and $i_\ell$ is the destination of the walk.
It is worth recalling that
$(A^\ell)_{ij}$ is the number of distinct walks from $i$ to $j$ 
with length $\ell$.

\medskip
Throughout the paper, $I_n$ denotes the $n \times n$ identity matrix, but we omit the subscript 
and simply write $I$ when the size is clear from the context.
We denote $e_i$ the $i$th column of the identity matrix 
and $\uno$ the all-ones vector
of appropriate size.
For any $A\in\R^{n\times n}$ let $\rho(A)$ be
its spectral radius, and $\vec(A)$ denotes the vectorization of $A$, that is, the $(n^2)$-vector obtained by stacking the columns of $A$ one on top of the other. Jointly with the vectorization operator, the Kronecker product $\otimes$
is a familiar tool to formally reducing the solution of matrix equations
to that of linear systems. In fact, we recall the following identities
that will often be needed in what follows:
\begin{itemize}
\item $\vec(ABC) = (C^T\otimes A)\vec(B)$,
\item $(A\otimes B)(C\otimes D) = (AC \otimes BD)$.
\end{itemize}
Finally, in the sequel we will use a couple of times 
the following classical results from
Perron-Frobenius theory,
see e.g., \cite{HorJoh}
or \cite[chapter 8]{MeyerBook}.

\medskip
\begin{theorem}   \label{thm:PF}
Let $A\in\R^{n\times n}$ be 
a matrix with nonnegative entries.
If there exists a positive vector $v\in\R^n$ such 
that $Av = \lambda v$ for some scalar $\lambda > 0$ 
then $\rho(A) = \lambda$.
In addition, if $A$ is also irreducible
then $\rho(A)$ is a simple eigenvalue, and 
an associated eigenvector can be chosen with 
all positive entries.
\end{theorem}

% ----------------------------------------------
\subsection{Quantifying node similarity}   \label{sec:3}

Given a (directed) graph $\G = (V,E)$,
a node similarity measure 
is a symmetric function $\sigma:V\times V\mapsto \R$
that
quantifies how much any two nodes $i$ and $j$ are similar to each other, in some sense. The value $\sigma(i,j)$ is conveniently stored in the $(i, j)$ entry of an associated matrix, henceforth denoted $S$
and called similarity matrix.
Various node similarity measures have been introduced in the literature,
see e.g., \cite{Blondel+,LHN05,Pan+10} and
the ample reviews included in the PhD theses of 
Arnaud Browet \cite{BrowetPhD} and Thomas Cason \cite{CasonPhD},
but they 
are essentially derived by quantifying one of two kinds of node equivalence criteria: structural equivalence and regular equivalence \cite{HanRid,WasFau}.
Similarity measures based on structural equivalence 
are based in the immediate neighborhood of a node
and yield explicit formulas for the similarity matrix. 
Formally speaking, two nodes are structurally equivalent
if they have the same ascendants and descendants.
Clearly, this definition is too stringent in practice, since it entails perfect interchangeability of nodes. 
For example, no pair of nodes in the network 
\eqref{eq:6nodes} is structurally equivalent.
Therefore, 
this criterion is used in an analogical way, 
that is, considering the extent to which two given nodes satisfy the original definition based on their relationships with other nodes in the network.
For example, a similarity 
measure for two nodes $i$ and $j$ in a directed network 
can be defined by counting the number of parent and child nodes common to both $i$ and $j$: $\sigma(i,j) = |\Gamma^{in}_i \cap\Gamma^{in}_j| + |\Gamma^{out}_i \cap\Gamma^{out}_j|$.
In this case, the similarity matrix 
admits the simple explicit formula 
\begin{equation}   \label{eq:sim0BVD}
   S = AA^T + A^TA .
\end{equation}
This definition is not suitable when node degrees are highly variable, since it can take large values for nodes with high degrees even if only a small fraction of their neighbors coincide. 
A more appropriate measure for the case where node degrees are highly heterogeneous is one where the importance of parent and child nodes shared by two nodes is related to their incoming and outgoing degrees, e.g., as in the formula below: 
$$
   \sigma(i,j) = \frac{|\Gamma^{in}_i \cap\Gamma^{in}_j|}{d^{in}_id^{in}_j} + \frac{|\Gamma^{out}_i \cap\Gamma^{out}_j|}{d^{out}_id^{out}_j} ,
$$
where the fraction $0/0$ is interpreted as zero.
If all degrees are nonzero, the corresponding similarity matrix is given by
\begin{equation}   \label{eq:sim0F}
   S = D^{-1}_{out}AA^TD^{-1}_{out} + D^{-1}_{in}A^TAD^{-1}_{in} .
\end{equation}
Due to their local nature, similarity measures 
inspired by the structural equivalence criterion
may have significant shortcomings. For example, in a large graph two nodes may play similar roles even though they have no neighbors in common, and a similarity measure as the ones described above may produce inappropriate values \cite{LHN05}.
In fact, comparing the immediate neighborhood of two nodes 
%%% is often limited as it 
neglects the
full topology and complexity of the network. 
Such considerations lead to an extended definition of similarity called regular equivalence, where two nodes are considered to be similar if they are connected to other nodes that are themselves structurally or regularly similar.
For example, the partitioning of the nodes of the graph
\eqref{eq:6nodes}
into the sets $\{1,4\}$, $\{2,5\}$ and $\{3,6\}$
satisfies the regular equivalence principle.

\medskip 
Node similarities based on regular equivalence 
are recursive in nature and their computation may require considering walks 
with arbitrary length. For this reason, their formulas usually include a parameter which acts as a 
decay factor that dampens the contribution of longer walks. 
Similarity measures based on regular equivalence are shown in, e.g., \cite{Pan+10,JehWidom,Blondel+} for directed graphs and  in \cite{EstHig10,JehWidom,LHN05,Pan+10} for undirected graphs,
with applications to bibliometrics, 
graph matching, hypertext classification,
and recommender systems.
%%%  [for use in bibliometrics and hypertext classification (in effetti anche altri domini)]. 
%%% ---- Da Wikipedia ---
% Effectively, SimRank is a measure that says 
% "two objects are considered to be similar if 
% they are referenced by similar objects." 
%%% ----------------------
The computation of the associated similarity matrices usually requires
solving matrix equations, inverting matrices 
or computing matrix functions.
For example, 
the similarity matrix proposed in \cite{Blondel+}
is a nontrivial solution of the matrix equation $S = \gamma(ASA^T + A^TSA)$, where $\gamma>0$ is a parameter controlling the balancing of long walks relative to short ones. %%% [to be revised].
To understand how this works, consider that 
for any two nodes $i,j$ we have
$$
   S_{ij} = \gamma
   \bigg( \sum_{h: i\to h}\sum_{k:j\to k}  S_{jk}
   +  \sum_{h: h\to i}\sum_{k:k\to j}  S_{jk} 
   \bigg) .
$$
This identity defines
the similarity between $i$ and $j$ as being proportional to the total similarity of their respective parent and child nodes.
Also the SimRank algorithm computes similarities 
between pairs of
nodes within a directed graph based on the 
idea that ``two objects are similar if they are related to similar objects''  \cite{JehWidom}.
The corresponding similarity matrix is a solution of 
the matrix equation $S = c A^TSA$
where $c$ is a constant between 0 and 1.
Obviously, the computation of this type of measure has a high computational cost due to of the matrix algebra operations (products and inversions) required. For this reason, some authors have introduced similarity measures that limit the exploration of walks only up to a certain length, 
allowing the user to control the growth of the computational cost, see e.g., \cite{Bar+13}.

% ----------------------------------------------
\section{Generalized walk patterns}   \label{sec:GWP}

In this section, we revise the idea put forward by 
A. Browet and P. Van Dooren in \cite{BroVDo}
for solving role extraction problems in directed networks. 
To this purpose, we introduce hereafter the notion of generalized walk pattern, which is 
not present in \cite{BroVDo}
but may be useful to gain a different standpoint on that paper and to better understand the developments shown in the sequel.

\medskip
A walk pattern $\Psi$ is a 
nonempty sequence of letters $d$ and $r$, 
e.g., $\Psi = ddrrd$. We denote $|\Psi|$ the length
of the sequence, e.g., $|ddrrd| = 5$.
Let $\G = (V,E)$ be a digraph, and
let $\Psi = \psi_1\cdots\psi_\ell$ be a walk pattern, 
$\psi_k\in\{d,r\}$ for $k = 1,\ldots,\ell$.
A sequence of nodes $i_0,i_1,\ldots,i_\ell$
such that 
\begin{itemize}
\item
$(i_{k-1},i_k)\in E$ if $\psi_k = d$
\item
$(i_k,i_{k-1})\in E$ if $\psi_k = r$,
\end{itemize}
is called generalized walk with pattern $\Psi$, or $\Psi$-walk.
So the letters $d$ and $r$ indicate whether a link is 
traversed in the direct or reverse direction, respectively.
For example, an ordinary walk is a 
generalized walk with an all-$d$ pattern.
We say that $i_0,i_1,\ldots,i_\ell$ is 
a $\Psi$-walk from $i_0$ to $i_\ell$.
A walk pattern induces a matrix operator defined as follows.

\medskip
\begin{definition}   \label{def:1}
Let $\Psi$ be a walk pattern, 
$\Psi = \psi_1\cdots\psi_\ell$. 
Given two matrices $M,N\in\R^{n\times n}$ 
define the matrix
$\Psi(M,N) = X_{1}X_{2}\cdots X_\ell$
where $X_k = M$ if $\psi_k = d$ and
$X_k = N$ if $\psi_k = r$.
When $N = M^T$ we use the 
shorter notation $\Psi(M)$ in place of $\Psi(M,M^T)$.
\end{definition}

\medskip
Looking at the matrix operator associated with a 
walk pattern, it is not difficult to derive 
an algebraic expression 
for counting the generalized walks between each 
pair of nodes in a digraph.

\medskip
\begin{theorem}   \label{thm:Psi(A)}
Let $\Psi$ be a walk pattern. If $A$ is the adjacency matrix
of the digraph $\G$ then 
the $(i,j)$-entry of $\Psi(A)$ is the number of $\Psi$-walks from $i$ to $j$.
\end{theorem}

\begin{proof}
Proceed by induction on $|\Psi|$. 
If $|\Psi| = 1$ then either $\Psi = d$ 
and $\Psi(A) = A$, or
$\Psi = r$ and $\Psi(A) = A^T$. 
In both cases the claim is trivial, 
as a generalized walk of length $1$ reduces to a single arc. So let $\Psi = \psi_1\cdots\psi_n$
with $n>1$. 
Assume $\psi_n = d$ for definiteness, the other case being analogous. 
Thus $\Psi(A) = \Psi'(A)A$ where $\Psi' = \psi_1\cdots\psi_{n-1}$. By Definition \ref{def:1},
$$
   \Psi(A)_{ij} = \sum_{k=1}^n
   \Psi'(A)_{ik}A_{kj} 
   = \sum_{k:k\to j} \Psi'(A)_{ik}.
$$
By induction, $\Psi'(A)_{ik}$ is 
the number of $\Psi'$-walks from $i$ to $k$.
So the rightmost sum counts the number of 
sequences $i = i_0,i_1,\ldots,i_{n-1},i_n = j$
where $i_{n-1} = k$ and $k\to j$,
that is, the $\Psi$-walks from $i$ to $j$.
\end{proof}

% ----------------------------------------------

\subsection{The Neighbourhood Pattern Similarity}
\label{sec:BVD}

From here on, the symbol $A$ denotes the
adjacency matrix of a digraph $\G = (V,E)$.
For any fixed walk pattern $\Psi$ and $i,j=1,\ldots,n$
let $s^\Psi_{ij}$ be the number 
defined as follows:
\begin{equation}   \label{eq:kPsiij}
   s^\Psi_{ij} = \sum_{k=1}^n [\Psi(A)]_{ik}[\Psi(A)]_{jk}
   = [\Psi(A)\Psi(A)^T]_{ij} .
\end{equation}
This number counts the number of common destinations
of $\Psi$-walks starting from $i$ and $j$, 
where each destination is repeated as many times as there are different ways to reach it.
The Neighbourhood Pattern Similarity of nodes $i$ and $j$
is defined in \cite{BroVDo,Marchand+} as the number
\begin{equation}   \label{eq:NPS}
   \sigma(i,j) = \sum_{\ell = 1}^\infty \sum_{|\Psi| = \ell} 
   \beta^{2\ell-2} 
   s^\Psi_{ij} ,
\end{equation}
where $0 \leq \beta < 1$ is a parameter.
The role of $\beta$ is twofold, namely, 
to lessen the relevance of longer $\Psi$-walks and
to ensure the convergence of the series.
In fact, we will see shortly that if $\beta$ is sufficiently small then the series \eqref{eq:NPS}
converges for all $i,j=1,\ldots,n$.
It is not difficult to recognize from \eqref{eq:kPsiij} 
that $\sigma(i,j)$ is the $(i,j)$-th entry of the matrix 
\begin{equation}   \label{eq:BVDmatrix}
   S = \sum_{\ell = 1}^\infty \sum_{|\Psi| = \ell} 
   \beta^{2\ell-2} \Psi(A)\Psi(A)^T .
\end{equation}
This matrix is actually a similarity matrix based on 
a regular equivalence criterion that extends 
the similarity matrix defined in \eqref{eq:sim0BVD}
beyond the immediate node neighborhoods.
Indeed, when $\beta = 0$ then the 
series \eqref{eq:BVDmatrix} reduces to $AA^T + A^TA$,
that is, the matrix in \eqref{eq:sim0BVD}.
The intuition behind \eqref{eq:BVDmatrix} is 
that a pair of nodes is highly similar 
if they can reach many common destinations by 
generalized walks of the same kind and length.
Informally put, while $(AA^T + A^TA)_{ij}$
counts the number of parent and child nodes common to both $i$ and $j$, $\sigma(i,j)$ is a weighted sum of 
not only the common parents and children, but also grandparents, uncles, nieces, etc.\ to both $i$ and $j$.
Additionally, the contribution of each common relative is weighted according to the number and length of the generalized walks that include it.
The following statement collects from \cite{BroVDo,Marchand+}
the main facts concerning the existence and computation 
of \eqref{eq:BVDmatrix}.

\medskip
\begin{theorem}   \label{thm:BVDmain}
Let $S_1 = AA^T+A^TA$ and
\begin{equation}   \label{eq:BVDfixpoint}
   S_{k+1} = S_1 + \beta^2 (A S_k A^T + A^T S_k A)
\end{equation}
for $k \geq 1$. Then
$$
   S_k = \sum_{\ell = 1}^k \sum_{|\Psi| = \ell} 
   \beta^{2\ell-2} \Psi(A)\Psi(A)^T .
$$
The iteration \eqref{eq:BVDfixpoint}
is convergent if and only if
$$
    \beta^2 < \frac{1}{\rho(A\otimes A + A^T\otimes A^T)} .
$$
In this case, the limit $S = \lim_{k\to\infty} S_k$
is the sum of \eqref{eq:BVDmatrix} and is also 
the solution of the matrix equation
\begin{equation}   \label{eq:BVDmeq}
   S - \beta^2 (ASA^T + A^TSA) = S_1.
\end{equation}   
\end{theorem}

% ---------------------------------------------
\section{Main results}   \label{sec:main}

In Section \ref{sec:GWP} we considered a generalized 
version of the usual walks on a digraph allowing to traverse arcs both in the proper and in the opposite direction.
In this section we extend in a similar way the concept of 
simple random walk, introduce the counterpart of the Neighborhood Pattern Similarity that results from this extension, and
present the main results of this work.

\medskip
From here on, we suppose that the digraph $\G = (V,E)$
is strongly connected, that is, for any two nodes $i,j\in V$ 
there exists a walk from $i$ to $j$. Consequently, all 
in-degrees and out-degrees are positive, 
the diagonal degree matrices $D_{in}$ and $D_{out}$ are nonsingular and the adjacency matrix $A$ is irreducible.
Note that this assumption is not very restrictive, as it can be met by providing all nodes with a loop, i.e. a link from a node to itself, or by adopting one of the ways used to properly define the PageRank vector in an arbitrary network, see e.g., \cite{PageRank}.
Furthermore, in what follows we can also suppose that 
the digraph $\G$ 
is weighted and $A_{ij}$ is the (positive) weight
associated with the link $i\to j$. 
So, let $P = D_{out}^{-1}A$ and $Q = D_{in}^{-1}A^T$.
The matrices $P$ and $Q$ are
the (row stochastic) transition matrices of the random walk 
on $\G$ and on the graph obtained from $\G$ by reversing its links, respectively. 
Both $P$ and $Q$ are irreducible, by the irreducibility of $A$.

\medskip
Now, consider the matrix sequence
$\{S^{(\ell)}, \ell = 1,2\ldots\}$ defined as follows:
%%% Let $S^{(1)} = I$ and, for $\ell \geq 1$,
\begin{equation}   \label{eq:def_Sell}
   S^{(\ell)} = 2^{-\ell}\sum_{|\Psi| = \ell} 
   \Psi(P,Q) \Psi(P,Q)^T .
\end{equation}
We defer to \S\ref{sec:rw} the interpretation of this matrix.

\medskip
\begin{lemma}   \label{lem:4}
For $\ell\geq 1$ the matrix in \eqref{eq:def_Sell} 
verifies the identity 
\begin{equation}   \label{eq:rec_Sell}
   S^{(\ell)} = \frac12 \big( PS^{(\ell-1)}P^T
   + QS^{(\ell-1)}Q^T \big) 
\end{equation}
where we set $S^{(0)} = I$.
\end{lemma}

\begin{proof}
We proceed by induction on $\ell$.
If $\ell = 1$ then either $\Psi = d$ or $\Psi = r$.
Hence, by Definition \ref{def:1},
\begin{align*}
   S^{(1)} & = 
   \frac12 \sum_{|\Psi| = 1} \Psi(P,Q)\Psi(P,Q)^T  \\
   & = \frac12 \big(d(P,Q)\, d(P,Q)^T + r(P,Q)\, r(P,Q)^T \big)
   = \frac12 \big( PP^T + QQ^T \big) ,
\end{align*}  
which is \eqref{eq:rec_Sell} as $S^{(0)} = I$ by hypothesis. 
Now consider the $\ell > 1$ case. Any walk pattern $\Psi$ with $|\Psi| = \ell > 1$
has the form $\Psi = \psi_1 \Psi'$ where $\psi_1\in\{d,r\}$
and $\Psi'$ is a walk pattern with $|\Psi'| = \ell - 1$. 
If $\psi_1 = d$ then $\Psi(P,Q) = P\Psi'(P,Q)$,
otherwise $\Psi(P,Q) = Q\Psi'(P,Q)$.
By an inductive argument, from \eqref{eq:def_Sell} we have
\begin{align*}
   S^{(\ell)}    
   & = 2^{-\ell} \bigg[\sum_{|\Psi| = \ell, \psi_1 = d}
   \Psi(P,Q)\Psi(P,Q)^T + \sum_{|\Psi| = \ell, \psi_1 = r}
   \Psi(P,Q)\Psi(P,Q)^T \Bigg] \\
   & = 2^{-\ell} \bigg[
   \sum_{|\Psi| = \ell-1} P\Psi(P,Q)\Psi(P,Q)^T P^T + 
   \sum_{|\Psi| = \ell-1} Q\Psi(P,Q)\Psi(P,Q)^T Q^T   
   \Bigg] \\
   & = \frac12 \big( 
   PS^{(\ell-1)}P^T + QS^{(\ell-1)}Q^T \big) ,
\end{align*}
and the proof is complete.
\end{proof}

In what follows, we analyze
existence and computation of the matrix
\begin{equation}   \label{eq:defS*}
   S^* = \sum_{\ell = 1}^\infty \beta^{2\ell-2} S^{(\ell)} ,
\end{equation} 
which is precisely the matrix that 
we are going to introduce as similarity matrix.
Together with \eqref{eq:def_Sell}, this definition 
is apparently related to \eqref{eq:BVDmatrix}.
Hereafter, we prove that the partial sums 
of \eqref{eq:defS*} verify a recurrence
that mirrors the one in \eqref{eq:BVDfixpoint}.

\begin{lemma}   \label{lem:3}
Let $S_1 = PP^T+QQ^T$. 
For $k \geq 2$ the partial sum 
$S_k = \sum_{\ell=1}^k \beta^{2\ell-2} S^{(\ell)}$
verifies the recurrence
\begin{equation}   \label{eq:rec_Sk}
   S_k = S_1 +\frac{\beta^2}{2}
   (PS_{k-1}P^T + QS_{k-1}Q^T) .
\end{equation}
\end{lemma}
   
\begin{proof} 
Using \eqref{eq:rec_Sell} we obtain
\begin{align*}
   S_{k} = S_1 + \sum_{j=2}^{k} \beta^{2j-2} S^{(j)} 
   & = S_1 + \frac{1}{2} \sum_{j=2}^{k} \beta^{2j-2}
   \big( PS^{(j-1)}P^T + QS^{(j-1)}Q^T \big) \\
   & = S_1 + \frac{\beta^2}{2} \bigg[ 
   P \bigg[ \sum_{j=1}^{k-1} \beta^{2j-2} S^{(j)} 
   \bigg] P^T +
   Q \bigg[ \sum_{j=1}^{k-1} \beta^{2j-2} 
   S^{(j)} \bigg] Q^T \bigg] \\
   & = S_1 + \frac{\beta^2}{2} \big( 
   P S_{k-1} P^T + Q S_{k-1} Q^T \big) ,
\end{align*}
and the proof is over.
\end{proof}

The next theorem gives the conditions for the existence of the matrix \eqref{eq:defS*}
and characterises it as the solution of a
matrix equation, analogously to Theorem \ref{thm:BVDmain} for the matrix
\eqref{eq:BVDmatrix}.

\medskip
\begin{theorem}   \label{thm:main}
The series \eqref{eq:defS*} is convergent
if and only if $\beta^2 < 1$. 
In this case, $S^*$ is a symmetric, positive semidefinite matrix with positive entries, and is the unique solution of the matrix equation 
\begin{equation}   \label{eq:main}
   S - \frac{\beta^2}{2}(PSP^T + QSQ^T) = S_1 ,
\end{equation}
where $S_1 = PP^T + QQ^T$.
\end{theorem}

\begin{proof}
Writing the recurrence 
\eqref{eq:rec_Sk}
in vectorized form, we obtain
\begin{equation}   \label{eq:vecSk}
   \vec(S_k) = \vec(S_1) + \frac{\beta^2}{2}
   (P\otimes P + Q \otimes Q) \vec(S_{k-1}) .
\end{equation}
Note that convergence of the series \eqref{eq:defS*}
corresponds to convergence of this iteration.
Let $M = (P\otimes P + Q \otimes Q)/2$.
By well-known results in the analysis of 
linear fixed-point iterations,
the iteration
\eqref{eq:vecSk} is convergent for all initial vectors
if and only if $\beta^2\rho(M) < 1$.
It is not difficult to check that $\rho(M) = 1$.
Indeed, both $P\otimes P$ and $Q \otimes Q$ are 
irreducible row-stochastic matrices, since both $P$ and $Q$ are.
So $M$ is irreducible and row stochastic. Indeed,
let $\uno = (1,\ldots,1)^T$. From $P\uno = Q\uno = \uno$
we derive $(P\otimes P)(\uno\otimes \uno) = 
(Q\otimes Q)(\uno\otimes \uno) = \uno\otimes \uno$, and
the identity
$$
   \frac12(P\otimes P + Q \otimes Q)(\uno\otimes \uno) = 
   \uno\otimes \uno
$$
proves that $1$ is an eigenvalue of $M$
associated to a positive eigenvector.
Thus $\rho(M) = 1$ follows from Theorem \ref{thm:PF}. Consequently, if
$\beta^2 < 1$ then \eqref{eq:vecSk} is convergent.
To prove the converse implication 
it is sufficient to show that 
$\vec(S_1)$ has a nonzero component along the 
dominant eigenspace of $M$, owing to the linearity of 
the recurrence \eqref{eq:vecSk}. 
So, let $v\in\R^{(n^2)}$ be a left eigenvector of 
$M$ associated with 
its spectral radius. By Theorem \ref{thm:PF}
we can safely assume that 
all entries of $v$ are positive since $M$ is irreducible.
Hence $v^T\vec(S_1) > 0$
and we are done.

\medskip
We address now the second part of the claim.
Symmetry and positive definiteness of $S^*$ follow from
\eqref{eq:defS*}, as all terms in the series are 
symmetric and positive semidefinite by virtue of Lemma \ref{lem:4}. By hypothesis, for every $i,j\in V$ there exists al least
one generalized walk from $i$ to $j$, hence $S^{(k)}_{ij} > 0$ for some $k>0$.
Thus all entries of $S^*$ are positive.
The identity $S^* - (\beta^2/2)(PS^*P^T + QS^*Q^T) = S_1$
follows from \eqref{eq:rec_Sk} by a continuity argument. 
To prove that $S^*$ is characterized by this identity
we look at the vectorized form of it, namely,
$$
   (I_m - \beta^2 M) \vec(S^*) = \vec(S_1) ,
$$
where $m = n^2$ and $M = (P\otimes P + Q \otimes Q)/2$ as before.
If $\beta^2 < 1$ then all eigenvalues of $\beta^2 M$
lie strictly inside the unit circle in the complex plane,
as $\rho(M) = 1$. This implies that the matrix 
$I_m - \beta^2 M$ is invertible,
since its eigenvalues lie inside the 
ball with center in 1 and radius $\beta^2 \rho(M) < 1$,
which excludes zero.
Hence, the solution to \eqref{eq:main} is unique.
\end{proof}

% ---------------------------------------------
\subsection{Solving \eqref{eq:main} by iteration}

The next result provides a globally convergent iterative 
method for solving \eqref{eq:main}, together with a
convergence bound. 

\medskip
\begin{corollary}   \label{cor:solving}
Let $\beta^2 < 1$. The iteration
$$
   Z_{k} = P(I + (\beta^2/2) Z_{k-1})P^T +  
   Q(I + (\beta^2/2) Z_k)Q^T
$$
is globally convergent to the
matrix $S^*$ in \eqref{eq:defS*} and
$$
   \|Z_{k} - S^*\|_{\max} \leq 
   \beta^2 \|Z_{k-1} - S^*\|_{\max} ,
$$
where $\|\,\cdot\,\|_{\max}$ denotes the matrix max-norm,
$\|X\|_{\max} = \max_{i,j}|X_{ij}|$.
Moreover, if $\|Z_k - Z_{k-1}\|_{\max}\leq \varepsilon$
then $\|Z_k - S^*\|_{\max} \leq \varepsilon\beta^2/(1-\beta^2)$.
\end{corollary}

\begin{proof}
Let $E_k = Z_k - S^*$. 
Rewriting \eqref{eq:main} as 
$S = P(I-(\beta^2/2)S)P^T + Q(I-(\beta^2/2)S)Q^T$
we easily obtain 
$E_{k} = (\beta^2/2) (PE_{k-1}P^T + QE_{k-1}Q^T)$. In vectorized form,
$$
   \vec(E_{k}) = (\beta^2/2)(P\otimes P + Q \otimes Q)\vec(E_{k-1}) .
$$
Let $M = (P\otimes P + Q \otimes Q)/2$.
Thus $M$ is nonnegative and all the row sums are equal to $1$,
so $\|M\|_\infty = 1$. Taking norms in the identity above we have
$$
   \| \vec(E_{k})\|_\infty \leq \beta^2\|M\|_\infty 
   \| \vec(E_{k-1})\|_\infty = \beta^2 
   \| \vec(E_{k-1})\|_\infty .
$$
To complete the first part of the proof it is sufficient to observe that
$\| \vec(E_k)\|_\infty = \|E_k\|_{\max}$.
Finally,
\begin{align*}
   \|Z_k - S^*\|_{\max} & \leq
   \beta^2 \|Z_{k-1} - S^*\|_{\max} \\
   & \leq \beta^2 \big( \|Z_k - S^*\|_{\max} +
   \|Z_k - Z_{k-1} \|_{\max} \big) \\
   & \leq \beta^2 (\|Z_k - S^*\|_{\max} + \varepsilon) ,
\end{align*}
and the last claim follows by rearranging terms.
\end{proof}

The iteration in Corollary \ref{cor:solving}
is equivalent but somewhat more efficient than \eqref{eq:rec_Sk} 
in that it avoids managing the matrix $S_1 = PP^T+QQ^T$
in the right-hand side.
Algorithm \ref{alg:1}
provides a description of the overall computational procedure, endowing the iteration with a stopping criterion based on the stepsize.

\begin{algorithm}[t]   \label{alg:1}
  \caption{Iterative solution of \eqref{eq:main}}
  \KwIn{$A\in\R^{n\times n}$ (adjacency matrix), 
  $\beta\in(0,1)$, $\varepsilon> 0$ (stopping criterion)}
  \KwOut{$S$, approximate solution of \eqref{eq:main}}
  Compute $P = \mathrm{Diag}(A\uno)^{-1}A$ and 
  $Q = \mathrm{Diag}(A^T\uno)^{-1}A^T$\;
  $Z_{new} = O$\;
  \Repeat{$\|Z_{old} - Z_{new}\| \leq \varepsilon$}{
     $Z_{old} = Z_{new}$\;
     $Z_{new} = I + (\beta^2/2)Z_{new}$\;
     $Z_{new} = PZ_{new}P^T + QZ_{new}Q^T$\;   
  }
  $S = Z_{new}$\;     
\end{algorithm}

% ----------------------------------------------
\subsection{The limiting behavior}

In this paragraph we look more closely 
at the dependence on $\beta$ of the matrix \eqref{eq:defS*}.
So, define the matrix-valued function 
$S(\beta) = \sum_{k=1}^\infty \beta^{2k-2}S^{(k)}$. 
It is immediate to see that $S(0) = PP^T+QQ^T$,
which is the same matrix as in \eqref{eq:sim0F}.
Moreover,
Theorem \ref{thm:main} shows that $S(\beta)$ is well
defined if and only if
$\beta^2< 1$, thus it is interesting to see 
how $S(\beta)$ behaves when $\beta^2$ approaches $1$. 
This issue is resolved by the result here below.

\begin{theorem}
In the notations above,
$$
   \lim_{\beta^2\to 1^-} (1-\beta^2) S(\beta) = 
   2\,\mathrm{trace}(W) \uno\uno^T ,
$$
where $W$ is the solution of the homogeneous equation
$$
   \frac12 (P^TWP+Q^TWQ) = W
$$ 
scaled so that $\sum_{i,j}W_{ij} = 1$.
\end{theorem}

\begin{proof}
The matrix $M = (P\otimes P + Q\otimes Q)/2$ is irreducible and row stochastic. 
Hence, $M$ can be decomposed using the Jordan canonical form 
$$
   M = VJV^{-1} = V \begin{pmatrix}
   1 & 0 \\ 0 & \hat J \end{pmatrix} V^{-1} ,
$$
where $J$ is the Jordan matrix of $M$ with the Perron eigenvalue
placed in the top-leftmost corner
for notational convenience.
By appropriately scaling the matrix $V$, 
its first column $v_1$ is the all-ones vector in $\R^{(n^2)}$
and the first column $w_1$ of $V^{-T}$ is a left Perron eigenvector of $M$ scaled so that $w_1^Tv_1 = 1$.
%  (equivalently, the dominant right eigenvector of AT ).
Hence we have the decomposition $M = v_1w_1^T + M_0$ 
where $M_0$ is a matrix such that $\rho(M_0) = \rho(\hat J) \leq 1$
and $v_1w_1^TM_0 = M_0v_1w_1^T = O$, 
while the rank-one matrix $v_1w_1^T$ is idempotent.

\medskip
For $\beta^2 < 1$ we have 
$\vec(S(\beta)) = (I-\beta^2M)^{-1}\vec(S_1)$.
Considering the Neumann series expansion
$(I - \beta^2 M)^{-1} = \sum_{k=0}^\infty \beta^{2k}M^k$,
which is convergent for $\beta^2<1$,
the vectorization of $S(\beta)$ 
can be expanded in series of $\beta^2$ as follows:
\begin{align*}
   \vec(S(\beta)) = (I_{m} - \beta^2 M)^{-1}\vec(S_1) 
   & = \sum_{k=0}^\infty \beta^{2k} 
   \big[ v_1w_1^T + M_0\big]^k \vec(S_1) \\
   & = \sum_{k=0}^\infty \beta^{2k} \big[
   v_1w_1^T + M_0^k \big] \vec(S_1) \\
   & = \frac{w_1^T \vec(S_1)}{1-\beta^2} v_1 + 
   (I-\beta^2 M_0)^{-1} \vec(S_1) .
\end{align*}
Inverting the vec operator and 
multiplying by $1-\beta^2$ we get
$$
   (1-\beta^2)
   S(\beta) = w_1^T \vec(S_1)
   \uno\uno^T + (1-\beta^2) X(\beta),
$$
where $X(\beta)\in\R^{n\times n}$ is the matrix such that
$\vec(X(\beta)) = (I_m-\beta^2 M_0)^{-1} \vec(S_1)$.
However, $\lim_{\beta^2\to 1} (I_m-\beta^2 M_0)^{-1} \vec(S_1) = \vec(S_1)$, so
$$
   \lim_{\beta^2\to 1^-} (1-\beta^2) S(\beta) = 
   w_1^T \vec(S_1) \uno\uno^T .
$$
Furthermore, let $W = \vec^{-1}(w_1)$.
From the equations $M^Tw_1 = w_1$ and $w_1^Tv_1 = 1$
we conclude that $W$ is the solution of the matrix equation
$(P^TWP + Q^TWQ)/2 = W$ scaled so that 
$\uno^TW\uno = 1$. 
Finally, 
\begin{align*}
   w_1^T \vec(S_1) = \mathrm{trace}(WS_1)
   & = \mathrm{trace}(W(PP^T+QQ^T)) \\
   & = \mathrm{trace}(P^TWP+ Q^TWQ)
   = \mathrm{trace}(2W) ,
\end{align*}
and the proof is complete.
\end{proof}

The previous theorem reveals that
$S(\beta)$ becomes increasingly uniform as $\beta^2$ approaches 1, making it more difficult to accurately recover node roles
from its rows and columns. 
On the other hand, when $\beta^2$ is `small' then the
computed similarity scores approach those
obtained by the matrix in \eqref{eq:sim0F};
for comparison, recall that the analogous limit 
for the 
similarity matrix 
\eqref{eq:BVDmatrix}
is the matrix in \eqref{eq:sim0BVD}.

% -------------------------------------------------------
\subsection{Probabilistic interpretation}   \label{sec:rw}

The fact that the matrix $S^*$ in
\eqref{eq:defS*} is built defined in terms of stochastic matrices suggests that the resulting pairwise node similarities can be interpreted in terms of
Markov chains
\cite{Seneta}. In fact, consider the following definition.

\medskip
\begin{definition}
Let $\Psi$ be a walk pattern.
A generalized random $\Psi$-walk on $\G = (V,E)$ is a sequence of $|\Psi|$ nodes visited by a walker that starts at a given node and repeatedly moves to a randomly chosen node, so that the resulting node sequence 
is a generalized walk with pattern $\Psi$. 
At each step, the next vertex is chosen at random among the possible candidates 
in proportion to the weight of the (incoming or outgoing) arcs.
\end{definition}

\medskip
Thus a generalized random walk on $\G$ is a finite-length,
non-homogeneous
Markov chain where the transition matrices may alternate between $P$ and $Q$, depending on the given walk pattern. 
The next lemma is an immediate consequence of 
the definition above, see \cite{Seneta}.

\medskip
\begin{lemma}   \label{lem:rw}
Let $\Psi$ be a walk pattern.
Then the matrix $\Psi(P,Q)$ is a row stochastic matrix
and $\Psi(P,Q)_{ij}$ is the probability  
that a random walker initially placed in $i$
ends in $j$ 
after a random $\Psi$-walk.
\end{lemma}

\medskip
A probabilistic interpretation of the entries of $S^{(\ell)}$ in \eqref{eq:def_Sell} is shown below.

\begin{theorem}
The entry $S^{(\ell)}_{ij}$ is the probability that two 
generalized random walks 
with the same pattern
chosen uniformly at random among all patterns of length $\ell$
starting from $i$ and $j$ end up at the same node.
\end{theorem}

\begin{proof}
From Lemma \ref{lem:rw} and the identity
$$
   \big[\Psi(P,Q)\Psi(P,Q)^T\big]_{ij}
   = \sum_{k=1}^n \Psi(P,Q)_{ik}\Psi(P,Q)_{jk} ,
$$ 
it is not difficult to recognize that 
the $(i,j)$-entry of the matrix $\Psi(P,Q)\Psi(P,Q)^T$
is the probability that two walkers starting at $i$ and $j$
will meet at a common destination 
after performing a generalized random walk 
with pattern $\Psi$.
Given that there are exactly $2^\ell$ distinct
walk patterns with length $\ell$,
$S^{(\ell)}_{ij}$ is the expected value
of that probability
when $\Psi$ is chosen uniformly at random given its length
$\ell$.
\end{proof}

\medskip
Therefore, for any $\ell \geq 1$ 
the function $\sigma^{(\ell)}(i,j) = S^{(\ell)}_{ij}$
is a similarity measure based on the following procedure:
\begin{itemize}
\item
Two walkers are initially placed on nodes $i$ and $j$.
\item
They choose uniformly at random a walk pattern $\Psi$ 
with $|\Psi| = \ell$ and perform a 
generalized random $\Psi$-walk.
\item
The value of $\sigma^{(\ell)}(i,j)$ is the probability 
that they meet at the end of the walk.
\end{itemize}
Finally, the matrix $S^*$ in \eqref{eq:defS*}
defines the similarity 
$\sigma^*(i,j) = \sum_{\ell = 1}^\infty 
\beta^{2\ell-2}\sigma^{(\ell)}$
by summing up the contribution of 
generalized random walks with arbitrary length. 
In a nutshell, in passing from 
\eqref{eq:BVDmatrix} to \eqref{eq:defS*} we are replacing the 
enumeration of $\Psi$-walks with 
the probability of random $\Psi$-walks.

% ----------------------------------------------
\section{A special property of $S^*$}
\label{sec:SBM}

This section aims to prove an invariance property of the similarity matrix $S^*$ that demonstrates its 
robustness to significant variations in node degrees.
This property reveals itself when the 
stochastic matrices $P$ and $Q$ arise from
the average matrix of a random matrix ensemble
known as degree-corrected stochastic block model.
A brief introduction to this model is given in the next paragraph for a better understanding of the background.
The result referred to above is presented in the subsequent paragraph.

% ----------------------------------------------
\subsection{Stochastic block models}

The Stochastic Block Model (SBM) 
is one of the earliest and 
widespread techniques to describe and recognize 
roles and group structures in graphs and networks
\cite{SBMeasy,dc-SBM}.
It assumes that the nodes of a network can be 
partitioned into disjoint
groups or blocks, such that the probability of a link between two nodes depends only on the blocks to which they belong. 
Basically, an SBM with $r$ blocks 
for a network with $n$ nodes consists of 
a role matrix $B\in\R^{r\times r}$
with entries in $[0,1]$ and a membership matrix 
$\Theta\in\{0,1\}^{n\times r}$
that encodes the block structure.
The number $B_{pq}$ denotes the probability 
that there is a link from any node in block $p$
to any node in block $q$.
The $i$-th row of the matrix $\Theta$ 
contains exactly one nonzero entry, whose position indicates which block the $i$-th node belongs to.
The $(i,j)$-entry of the adjacency matrix 
of a network in this model
is assumed to be a Bernoulli random variable 
with expectation equal to the $(i,j)$-entry
of the matrix $\Theta B\Theta^T$.
Thus, SBMs can be used both as generative models, 
describing how a network with prescribed role structure could be randomly generated, and as a tool for understanding network structures by 
computing the model parameters $B$ and $\Theta$ 
that best describes a given adjacency matrix.   
SBMs are also widely used as a theoretical benchmark for 
the analysis of community detection algorithms. 
Concerning the neighborhood pattern similarity matrix
\eqref{eq:BVDmatrix}, 
the authors of \cite{Barbarino22} 
have shown that placing $A = \Theta B\Theta^T$
in the equation \eqref{eq:BVDmeq}
then the solution $S$ of that equation
permits the exact recovery of the blocks of the model,
under mild hypotheses.
This result refines a similar result found in 
\cite{Marchand+} for the case where the role matrix 
$B$ has binary entries, i.e.,
the block model is deterministic. Furthermore,
it is also shown in \cite{Barbarino22} 
that random adjacency matrices drawn from an SBM can provide reliable approximations to the blocks of the model
if they are sufficiently large.

\medskip
Stochastic block models as described above 
imply that all nodes within the same block have the same average degree. 
The degree-corrected SBM 
(dc-SBM) is a variant of the basic SBM designed to allow a greater degree heterogeneity \cite{dc-SBM,FT_SM18}. 
A dc-SBM is parametrized by a role matrix 
$B\in\R^{r\times r}$ and a membership matrix 
$\Theta\in\{0,1\}^{n\times r}$ as above, and
two diagonal matrices with positive diagonal entries, $D_1$ and $D_2$, provided that
the entries of the matrix $A = D_1 \Theta B\Theta^TD_2$ belong to $[0,1]$. The effect of $D_1$ and $D_2$ is to 
allow large degree variations for nodes 
belonging to the same block.
As before, $A_{ij}$ is the mean value of 
a Bernoulli random variable associated with the existence of the link $(i,j)$ in the network.
In fact, $A$ represents the expected adjacency matrix of the random networks in the model.

% --------------------------------------
\subsection{The similarity matrix of the average dc-SBM matrix}   \label{sec:dcSBM}

In the following we prove that the similarity matrix
$S^*$ in \eqref{eq:defS*}
computed from the average matrix
$A = D_1\Theta B\Theta^T D_2$ of a dc-SBM 
can accurately recover the role structure of the model, 
independently of the degree correction 
introduced by $D_1$ and $D_2$.
Here, $B\in\R^{r\times r}$ is the role matrix
and
$\Theta\in\R^{n\times r}$ is the membership matrix.
In line with the 
strong connectedness assumption
in Section \ref{sec:main}, we require that no rows and columns of $B$ vanish.
To simplify the calculations that follow, we normalize the columns of $\Theta$
so that they are orthonormal vectors.
More precisely,
for $i = 1,\ldots r$ let $\chi^{(i)}\in\R^n$ be the characteristic vector of the $i$th block, that is,
$\chi^{(i)}_j = 1$ if $j$ is an index
in the $i$th block and $\chi^{(i)}_j = 0$ otherwise.
Moreover, let $n_i$ be the number of indices in the $i$th block,
so $n_1 + \cdots+n_r = n$.
Then the $i$th column of $\Theta$ is
\begin{equation}   \label{eq:Bcolumns}
   \Theta e_i = \frac{1}{\sqrt{n_i}} \chi^{(i)} , \qquad
   i = 1,\ldots,r.
\end{equation}
Think of $A$ as the adjacency matrix of a weighted
digraph and consider the row-stochastic transition matrices
$P$ and $Q$ naturally associated with $A$,
as in Section \ref{sec:main}. In other words, 
$P$ and $Q$ are defined by scaling the rows of 
$A$ and $A^T$, respectively, so that they become row stochastic.
Then $P$ and $Q$ can be factorized as 
\begin{equation}   \label{eq:PQ}
   P = \wh D_1\Theta B \Theta^TD_2 , \qquad
   Q = \wh D_2\Theta B^T\Theta^TD_1 ,
\end{equation}
for some auxiliary diagonal matrices $\wh D_1, \wh D_2$
that are
characterized by the equations $P\uno = Q\uno = \uno$
where $\uno = (1,\ldots,1)^T\in\R^n$.
Optionally, we may suppose that the role matrix $B$
is invertible. This assumption entails that 
roles are clearly identifiable, because 
the connection pattern of each block is 
markedly different from the other blocks. 

\medskip
\begin{lemma}   \label{lem:X}
Let the matrices $A = D_1\Theta B \Theta^TD_2$, $P$ and $Q$ be as described above.
Consider the $r\times r$ matrices 
$M = \Theta^T P\Theta$ and $N = \Theta^T Q\Theta$.
If $\beta^2<1$ then the matrix equation
\begin{equation}   \label{eq:X}
   X - (\beta^2/2) (M X M^T + N X N^T ) = 
   \Theta^T(PP^T + QQ^T)\Theta 
\end{equation}
has a unique solution $X\in\R^{r\times r}$,
which is symmetric positive semidefinite.
Furthermore, if the role matrix $B$ is invertible then 
$X$ is positive definite.
\end{lemma}

\begin{proof}
Let $v = (\sqrt{n_1},\ldots,\sqrt{n_r})^T \in\R^r$.
Owing to \eqref{eq:Bcolumns} we 
have $\Theta v = \uno$ and $\Theta^T\uno = v$.
Hence,
$$
   Mv = \Theta^TP\Theta v = 
   \Theta^T P \uno = \Theta^T \uno = v .
$$
The identity $Nv = v$ follows by analogous passages.
Thus the vector $v$, which is entrywise positive,
is an eigenvector of both $M$ and $N$, which are
nonnegative. 
Theorem \ref{thm:PF} allows us to conclude that 
$\rho(M) = \rho(N) = 1$.
The existence, uniqueness and positive semidefiniteness of 
the solution $X$ of \eqref{eq:X}
can be obtained 
by reasoning as in the proof of Theorem \ref{thm:main}.
Finally, if $B$ is invertible then both $M$ and $N$ have full rank,
therefore both $MM^T$ and $NN^T$ are positive definite.
Thus also the right-hand side of \eqref{eq:X} is positive definite, from which the claim follows.
\end{proof}

\begin{theorem}   \label{thm:dcsbm}
Let $\beta^2<1$ and let $A$ be the average matrix of a 
degree-corrected stochastic block model, that is,
$A = D_1\Theta B \Theta^TD_2$ in the notation introduced above.
Then 
the matrix $S^*$ in \eqref{eq:defS*}
is $S^* = \Theta X \Theta^T$ 
where $X\in\R^{r\times r}$
is the solution of \eqref{eq:X}.
In particular, if $i$ and $j$ are two indices 
belonging to the same block then the corresponding 
rows and columns of $S^*$ are equal.
Moreover, if the role matrix $B$ is nonsingular then
$\mathrm{rank}(S^*) = k$. 
\end{theorem}

\begin{proof}
To prove the theorem it is sufficient to show that 
$S = \Theta X \Theta^T$ fulfills the identity
$$
    S - (\beta^2/2)(PSP^T + QSQ^T) = PP^T + QQ^T .
$$
To this aim we need to exploit a special structure 
of the row-stochastic matrices $P$ and $Q$. 
Considering $P$, from \eqref{eq:PQ} we have 
$$
   \uno = P\uno  = \widehat D_1 \Theta B \Theta^T D_2 \uno
   = \widehat D_1 u ,
$$
where $u = \Theta B \Theta^T D_2 \uno$ is a linear combination of the columns of $\Theta$, 
see \eqref{eq:Bcolumns}. Hence, if $i$ and $j$ are two indices belonging to the same block of the dc-SBM then 
$u_i = u_j$. Therefore, also the corresponding diagonal entries 
of $\widehat D_1$ are equal:
$(\widehat D_1)_{ii} = 1/u_i = 1/u_j = (\widehat D_1)_{jj}$. 
Note that $u_i> 0$ for $i=1,\ldots,n$
since $\Theta B\Theta^TD_2$ has nonnegative entries and no null rows.
From the above facts we deduce that 
the column space of $P$ is included in the column space of $\Theta$.
With a similar line of reasoning, the same conclusion can be drawn for the matrix $Q$.
We arrive at the identities $P = \Theta\Theta^T P$
and $Q = \Theta\Theta^T Q$, 
since the columns of $\Theta$ are orthonormal
and $\Theta\Theta^T$ is the orthoprojector onto their span. 
Now, recall the definitions $M = \Theta^T P \Theta$
and $N = \Theta^T Q\Theta$. Multiplying 
the right-hand side of \eqref{eq:X} 
by $\Theta$ from the left and by $\Theta^T$ from the right
and simplifying, we obtain
\begin{align*}
   \Theta \big( X - (\beta^2/2) 
   (M X M^T + N X N^T )\big)\Theta^T 
   & = S^* -(\beta^2/2)(P\Theta X\Theta^T P^T
   +  Q\Theta X\Theta^T Q^T  )\\
   & = S^* -(\beta^2/2)(PS^* P^T +  QS^* Q^T  ) .
\end{align*}
On the other hand, using again \eqref{eq:X}
and the identities $P = \Theta\Theta^T P$
and $Q = \Theta\Theta^T Q$, we also have
\begin{align*}
   \Theta \big( X - (\beta^2/2) 
   (M X M^T + N X N^T )\big)\Theta^T 
   & = \Theta \Theta^T (PP^T + QQ^T )\Theta\Theta^T \\
   & = PP^T + QQ^T  .
\end{align*}
So $S^* = \Theta X \Theta^T$ by Theorem \ref{thm:main}.
In particular,
if $i$ and $j$ are two indices 
belonging to the $k$th block for some $k\in\{1,\ldots,r\}$
then $\Theta^Te_i = \Theta^Te_j = e_k/\sqrt{n_k}$, so 
both $S^*e_i = S^*e_j$ and $e_i^TS^* = e_j^TS^*$ follow.
To complete the proof, suppose that $B$ is nonsingular.
Thus Lemma \ref{lem:X} guarantees that 
$X$ has full rank, and the identity 
$\mathrm{rank}(S^*) = r$ follows.
\end{proof}

It is worth noting that, according to the previous theorem, the matrix $S^* = \Theta X\Theta^T$ 
not only does not depend on the diagonal scalings
$D_1$ and $D_2$, but also has a rank exactly equal to the number of roles in the model.

% -------------------------------------------------------
\section{Numerical examples}   \label{sec:numerical}

In this section we show some examples of the use of the new similarity matrix to illustrate its application in role extraction problems, using both synthetic and real-world networks. 

% -------------------------------------------------------
\subsection{Random adjacency matrices from stochastic block models}

\begin{figure}   
\begin{center}
\includegraphics[width=15cm]{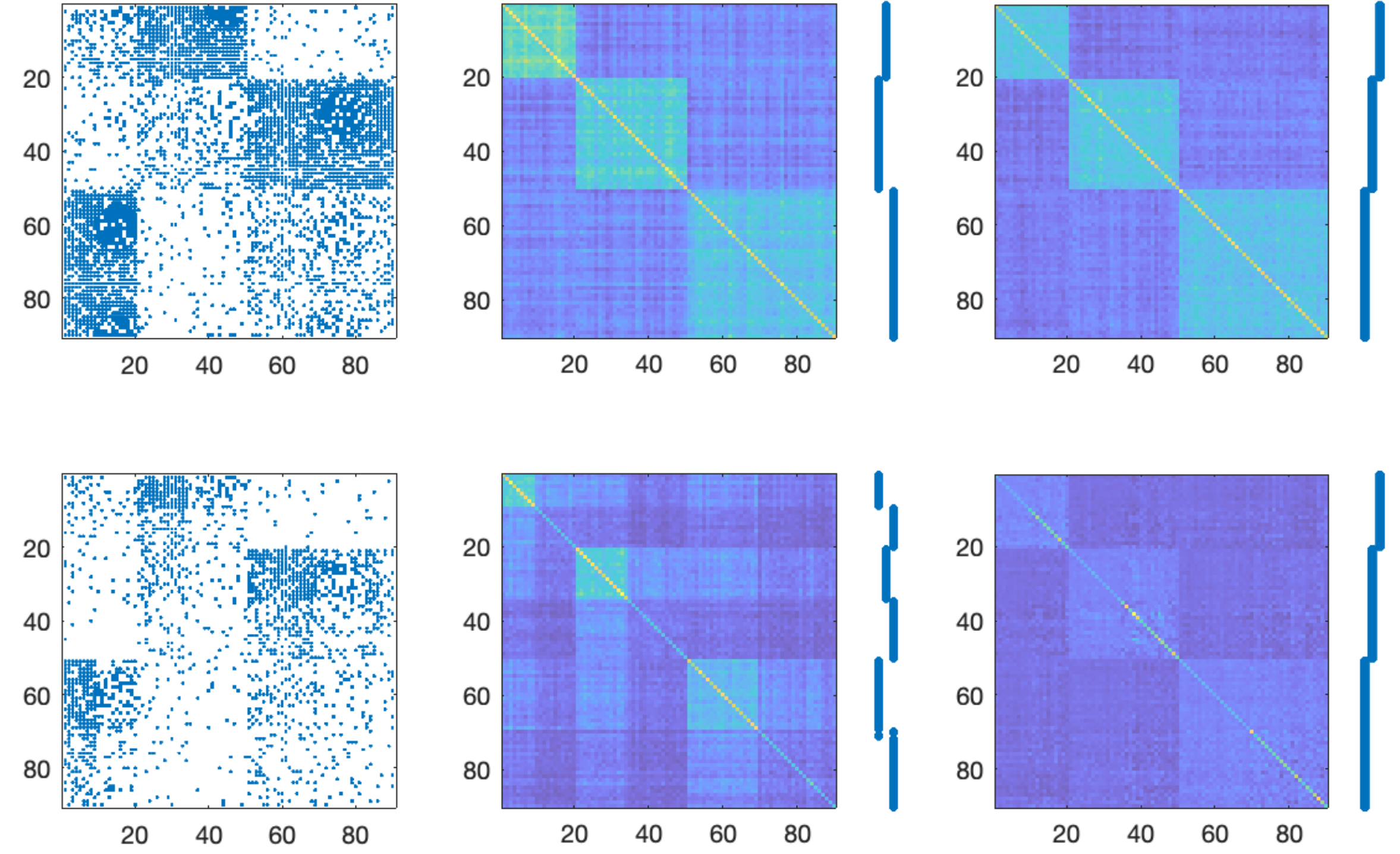}
\caption{Top row: A random adjacency matrix 
drawn from a Stochastic Block Model
(left) and the corresponding node similarity matrices $S$ (center) and $S^*$ (right), together with the cluster assignments 
computed by row clustering.
Bottom row: Same as above, but using a degree-corrected SBM. \label{fig:1}}
\end{center}
\end{figure}

An unsymmetric matrix $A$ was generated 
as a random adjacency matrix 
belonging to a stochastic block model with $3$ blocks.
The role matrix is
$$
   B = \begin{pmatrix}
   0.3 & 0.8 & 0.1 \\ 0.1 & 0.3 & 0.8 \\ 0.8 & 0.1 & 0.3
   \end{pmatrix} ,
$$
and the block sizes are $20$, $30$ and $40$.
Thus, the entry $A_{ij}$ is set to $1$
with probability $B_{pq}$ where $p$ is the block number
of $i$ and $q$ the one of $j$, otherwise $A_{ij} = 0$.
For a better visual effect, 
row and column indices are assigned to the blocks in sequence.
The top row of Figure \ref{fig:1} represents the 
matrix $A$ (left) and the similarity matrices $S$ from 
\eqref{eq:BVDmatrix} (center) and $S^*$ from
\eqref{eq:defS*} (right)
obtained from $A$.
The rows of these matrices were input to the 
$k$-means algorithm as implemented in the Matlab's function 
{\tt kmeans} with $k = 3$.
The computed clusters are indicated by the dot columns
on the side of each matrix, one dot per row.
In this example the connections among the three 
node groups are fairly uniform and the roles are perfectly detected from both $S$ and $S^*$.

\medskip
Next, another random adjacency matrix was generated 
from a degree-corrected version of the preceding block model.
More precisely, the role matrix and block 
assignments are as in the preceding example.
To increase degree heterogeneity,
the probability that a node receives or send a link 
is halved, for half the nodes in each block.  
The bottom row of Figure \ref{fig:1} shows 
the adjacency matrix (left) and the similarity matrices 
$S$ (center) and $S^*$ (right)
computed from it.
While the three roles can still 
be identified from the rows of $S^*$, 
the rows of $S$
lose their coherence and the $k$-means clustering 
applied to them no longer yield accurate results.

% -------------------------------------------------------
\subsection{The average dc-SBM matrix}

\begin{figure}   
\begin{center}
\includegraphics[width=15cm]{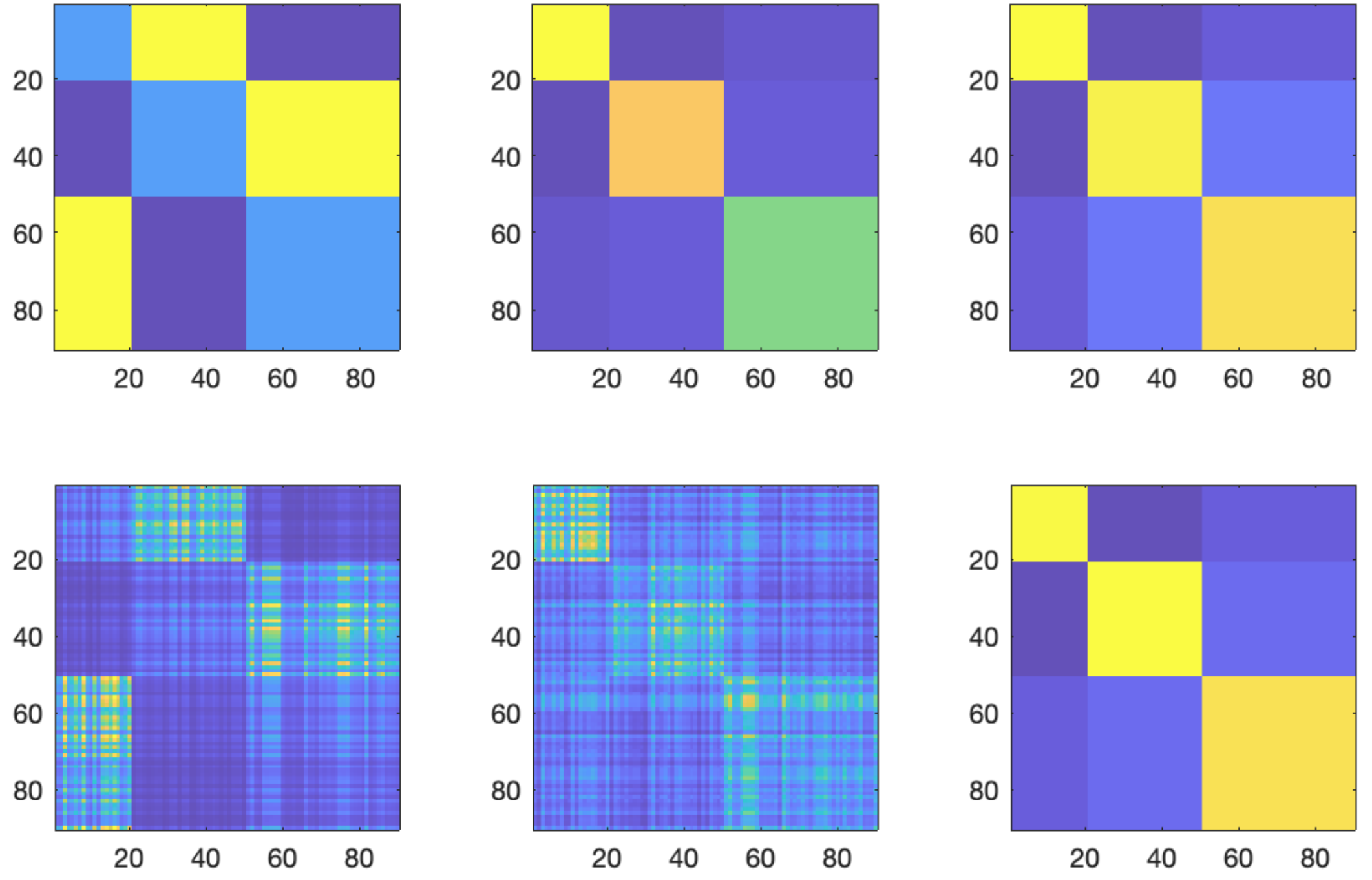}
\caption{Top row: The average matrix of a SBM
(left) and the corresponding similarity matrices $S$ (center) and $S^*$ (right). Darker colors represent smaller entries.
Bottom row: Same as above but with a degree-corrected SBM.\label{fig:2}}
\end{center}
\end{figure}

The next example is meant to illustrate 
the content of Section \ref{sec:SBM}
and, in particular, Theorem \ref{thm:dcsbm}.
Firstly, we consider a block partitioned
data matrix $A$ with constant blocks. This matrix is the average matrix of 
the SBM considered in the previous example,
that is, the values of the entries in each block 
are equal to the entry of $B$ corresponding to that block.
Looking at the top row of Figure \ref{fig:2},
the matrix $A$ is shown in colors on the leftmost panel,
along with the matrices $S$ (center) and $S^*$ 
(right) computed from $A$.
Both matrices $S$ and $S^*$ 
are block-partitioned with constant blocks, so that 
they can perfectly identify the blocks of the model. 
All matrices have rank 3.

\medskip
Then, we scaled the rows and column of $A$ by random coefficients, so that the resulting matrix is the average matrix of a degree-corrected SBM. 
As shown in the bottom row of Figure \ref{fig:2},
this operation degrades the matrix $S$, which becomes rank 6, while $S^*$ is unchanged, as expected from Theorem \ref{thm:dcsbm}.
The block structure of the matrix $S$ in the SBM case,
and
the rank doubling in the dc-SBM case, can be explained 
by the results in \cite{Barbarino22,Marchand+}.

% -------------------------------------------------------

\subsection{A faculty hiring network}

\begin{figure}   
\begin{center}
\includegraphics[width=15cm]{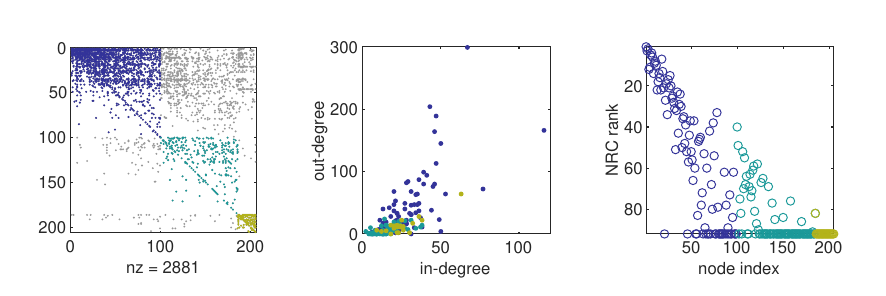}
\caption{Left: Adjacency matrix rearranged placing nodes
consecutively within each role. The diagonal blocks are in color to highlight position, size and consistency of each role.
Center: Scatter plot of node in-/out-degrees. 
Right: NRC rank vs node index. \label{fig:3}}
\end{center}
\end{figure}

The network analyzed in the next example is a directed and weighted network 
with 205 nodes and 2881 arcs
representing academic hiring relationships among 
universities in the US and Canada. Each node corresponds to a university or research institution and an arc from node $i$ to node $j$ indicates that faculty members at institution $j$ received their Ph.D.s from institution $i$. Arc weights indicate the number of people moving from $i$ to $j$.
This network is a subset of a larger dataset  encompassing nearly 19,000 tenure-track or tenured faculty placements across 461 North American academic units, collected between May 2011 and August 2013 as part of a large survey on 
institutional prestige, hierarchy and   
and gender issues in faculty hiring \cite{hiring-cs}.
The present analysis focuses specifically on the subset of the original data pertaining to Computer Science. 
Each institution in the network is assigned a
prestige score based on its ranking in a National Research Council (NRC) 
evaluation of doctoral and research programs available at the time of the survey.
The scores are available in the network meta-data
and vary from 1 to 92,
higher scores corresponding to lower prestige.

\begin{figure}   
\begin{center}
\includegraphics[width=15cm]{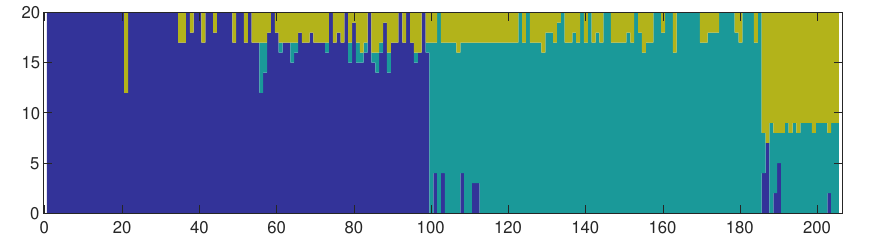}
\caption{Number of times each node in the network has been assigned a particular role by the $k$-means algorithm. The roles are represented by the same color scheme as in Figure \ref{fig:3}.   \label{fig:4}}
\end{center}
\end{figure}

\medskip
We clustered the rows of the similarity matrix $S^*$
using $k$-means, keeping the arc weights present 
in the original network data
and setting $\beta^2 = .2$.
The clustering was performed 20 times with $k = 3$ and
each node was assigned the role most frequently computed by the algorithm.
We renumbered the nodes within each role consecutively to improve visualization.
The left panel of Figure \ref{fig:3} shows the 
sparsity pattern of the adjacency matrix reordered in this way. 
The diagonal blocks of the matrix are colored in
blue, light blue and dark green
to make the position and size of each block more visible.
The same colors are used in the other panels of Fig. \ref{fig:3}.
The plot in the middle panel shows the 
distribution of node in-degrees and out-degrees,
illustrating the degree spread in the network.
The rightmost chart shows the NRC ranking of each node. It is clearly visible that the high ranked institutions belong to the first group, while the third group consists mainly of low ranked institutions.
The number of times each node has been assigned a particular role by the $k$-means algorithm is shown in Figure \ref{fig:4}.

\medskip
On the basis of the role partitioning thus obtained, 
we computed a plausible role matrix $B\in\R^{3\times 3}$
for a degree-corrected SBM that best describes the (unweighted) adjacency matrix of this network, 
based on the block structure shown in 
the left panel of Fig.\ \ref{fig:3}. The value of $B_{ij}$ 
is given by the ratio of nonzero entries in the 
block at position $(i,j)$ to the total number
of entries in the same block.
The matrix computed in this way is the following:
$$
   B = \begin{pmatrix}
    0.1508 & 0.0852 & 0.0763\\
    0.0095 & 0.0375 & 0.0052\\
    0.0131 & 0.0128 & 0.2800
    \end{pmatrix} .
$$
From the above observations, 
the nodes of this network seem to be organized in three groups with clearly defined roles. 
Top-level institutions belong to the first group, 
are selective in their recruitment and their graduates are widely disseminated throughout the education system.
A second group mainly consists of mid-level institutions, whose capacity to recruit staff 
is much less than that of the other groups
and whose graduates find little employment in the academic system. Finally, the universities in the last group are low ranked, somewhat attractive and their doctoral students are rarely recruited by institutions outside the same group.

% -------------------------------------------------------

\section{Conclusions}

We have described and analyzed an algorithm for
computing a node similarity matrix
for solving the role extraction problem
in directed networks.
This matrix is defined as the solution of a matrix equation,
for which we provide a convergent iteration,
and arises from a variation of the 
node similarity matrix originally introduced 
by A. Browet and P. Van Dooren in \cite{BroVDo} by replacing generalized $\Psi$-walks with generalized random $\Psi$-walks.
Theoretical analysis and numerical examples show that role extraction 
tasks based on the new similarity matrix achieve
remarkable performances also in networks with
heterogeneous node degree distributions. 
Furthermore, the definition of the new similarity matrix 
naturally extends to weighted networks.   

Admittedly, even if the adjacency matrix of the network is sparse, the new similarity matrix is completely full, as is 
normally the case with pairwise node similarity matrices 
derived from a regular equivalence criterion.
This makes the exact computation of such a matrix impractical for medium to large networks, due to the cubic computational cost of each single iteration of Algorithm \ref{alg:1}.
A similar problem concerning the calculation of 
the Browet-Van Dooren's similarity matrix 
%%% The analogous problem 
was solved %%% in \cite{BroVDo}
by endowing the corresponding iteration with a low-rank compression.
A similar technique could also be adopted here, e.g., by replacing the exact computation of the matrix $Z_{new}$
in Algorithm \ref{alg:1} with a low-rank approximation.
Other open questions in this work concern 
extending the definition of $S^*$ to graphs that have  nodes with no incoming or outgoing arcs, optimality of the number of roles and the parameter $\beta$.
Furthermore, the results in Section 
\ref{sec:dcSBM} may hold under weaker assumptions.
These developments could be the subject of further work.

% -------------------------------------------------------
\section*{Acknowledgements}

This work was supported in part by the 
Italian Ministry of University and Research
through the PRIN Project 20227PCCKZ
``Low Rank Structures and Numerical Methods in Matrix and Tensor Computations and their
Applications''.
The author is also affiliated to the INdAM-GNCS (Gruppo Nazionale di Calcolo Scientifico).
The author has no competing interests to declare that are relevant to the content of this article. 

\section*{Data Availability Statement}

The network data used to generate figures \ref{fig:3} 
and \ref{fig:4} are publicly available 
in the supplementary materials of the paper \cite{hiring-cs}.
%, see 
%{\tt https://www.science.org/doi/epdf/10.1126/sciadv.1400005}

\bibliographystyle{plain}
%%% \bibliography{roles}
% -----------------------------------------------

\medskip\noindent
\sc University of Udine, Udine, Italy.\\
\em Email address: {\tt dario.fasino@uniud.it}\\
\normalfont ORCID: 0000-0001-7682-0660
% -----------------------------------------------
\end{document}